\documentclass[12pt,thmsa]{article}
\usepackage{amsmath, latexsym, amsfonts, amssymb, amsthm, amscd}

\newtheorem{theorem}{Theorem}

\newtheorem{lemma}{Lemma}

\newcommand{\p}{\Bbb{P}}

\newcommand{\e}{\Bbb{E}}
\newcommand{\ex}{\Bbb{E}_x}

\newcommand{\ind}{\mbox{\rm 1\hspace{-0.04in}I}}

\newcommand{\R}{\mbox{\rm I\hspace{-0.02in}R}}

\newcommand{\ud}{\mathrm{d}}

\title{
\textbf{The Gapeev--K\"uhn stochastic game driven by a spectrally positive
L\'evy process.}}

\author{\textbf{E.J. Baurdoux\footnote{Department of Statistics, London School of Economics. Houghton street, {\sc London, WC2A 2AE, United Kingdom.} E-mail: e.j.baurdoux@lse.ac.uk}, A.E. Kyprianou\footnote{Department of Mathematical Sciences, University of Bath, Claverton Down, {\sc Bath, BA2 7AY, United Kingdom}. Email: a.kyprianou@bath.ac.uk, \, jcpm20@bath.ac.uk} }\footnote{Corresponding author.} \textbf{,\,J.C. Pardo$^{\dagger}$}}
\date{\footnotesize This version: \today}

\begin{document}

\maketitle

\begin{abstract}
\bigskip
In \cite{gak}, the stochastic game corresponding to perpetual convertible bonds   was considered when driven by a Brownian motion and a compound Poisson process with exponential jumps. We consider the same stochastic game but driven by a spectrally positive L\'evy process. We establish a complete solution to the game indicating four principle parameter regimes as well as characterizing the occurence of continuous and smooth fit. In \cite{gak}, the method of proof was mainly based on solving a free boundary value problem. In this paper, we instead use fluctuation theory and an auxiliary optimal stopping problem to find a solution to the game.
\\

\bigskip

\end{abstract}

\vspace{0.5cm}
\section{Introduction.}
Let $X=(X_t, t\geq 0)$ be a L\'evy process defined on a filtered probability space $(\Omega, \mathcal{F}, \mathbb{F}, \p)$, where $\mathbb{F}:=\{\mathcal{F}_t,t\ge 0\}$ is the filtration generated by $X$ which is naturally enlarged (see for instance Definition 1.3.38 in \cite{bic}).  For $x\in \R$ denote by $\p_x$ the law of $X$ when it is started at $x$ and write simply $\p_0=\p$. Accordingly, we shall write $\e_x$ and $\e$ for the associated expectation operators. In this paper we shall assume throughout that $X$ is \textit{spectrally positive} meaning here that it has no negative jumps and that it is not a subordinator. It is well known that the latter allows us to talk about the Laplace exponent $\psi(\theta):[0,\infty) \to \R$, i.e.
\[
\e\Big[e^{-\theta X_t}\Big]=:e^{t\psi(\theta)}, \qquad t, \theta\ge 0,
\]
and the Laplace exponent is given by the L\'evy-Khintchine formula
\begin{equation}\label{lk}
\psi(\theta)=\mu\theta+\frac{b^2}{2}\theta^2+\int_{(0,\infty)}\big(e^{-\theta x}-1+\theta x\ind_{\{x<1\}}\big)\Pi(\ud x),
\end{equation}
where $\mu\in \R$, $b^2\ge 0$ and $\Pi$ is a measure on $(0,\infty)$ called the L\'evy measure of $X$ and satisfies
\[
\int_{(0,\infty)}(1\land x^2)\Pi(\ud x)<\infty.
\]
The reader is referred to Bertoin \cite{Be} and Sato \cite{sa} for a complete introduction to the theory of L\'evy processes.

Denote by $\mathcal{T}_{0,\infty}$ the family of all $[0,\infty]$-valued stopping times with respect to $\mathbb{F}$. We are interested in establishing a solution to a special class of stochastic games which are driven by spectrally positive L\'evy processes.  Specifically, for  $ \alpha \ge 0$ and $\beta,q,K>0$, let
\[
L_t:=e^{-qt+X_t}+\int_0^t e^{-qs}(\alpha+\beta e^{X_s})\ud s,
\]
and
\[
U_t:=e^{-qt}(e^{X_t}\lor K)+\int_0^t e^{-qs}(\alpha+\beta e^{X_s})\ud s.
\]
We are interested in the stochastic game consisting of two players and expected pay-off given by
\begin{equation}\label{val}
M_x(\tau, \sigma):=\ex\Big[L_{\tau}\ind_{\{\tau<\sigma\}}+U_{\sigma}\ind_{\{\sigma\leq \tau\}}\Big],
\end{equation}
for $x\ge 0$. The inf-player's objective is to choose some $\sigma\in\mathcal{T}_{0,\infty}$ which minimizes (\ref{val}), whereas the sup-player chooses some $\tau\in\mathcal{T}_{0,\infty}$ which maximizes this quantity.  We are principally interested in showing  the existence of a stochastic saddle point (also known as Nash equilibrium cf. Ekstr\"om and Peskir \cite{EkP}). That is, we want to find $\tau^*$ and $\sigma^*$ such that
\begin{equation}\label{nash}
M_x(\tau,\sigma^*)\leq M_x(\tau^*, \sigma^*)\leq M_x(\tau^*, \sigma) \qquad \textrm{for all}\quad \tau, \sigma.
\end{equation}
Note that (\ref{nash}) implies for each $x$,
\begin{equation}\label{stgame}
V(x):=\sup_{\tau}\inf_{\sigma}M_x(\tau, \sigma)=\inf_{\sigma}\sup_{\tau}M_x(\tau, \sigma),
\end{equation}
the so called Stackelberg equilibrium. We call $V(x)$ the value of the game when it exists.
Note that for $x\geq \log K$
\[M_x(\tau,0)=e^x =M_x(0, 0)= M_x(0, \sigma)\]
for any $\tau,\sigma$, i.e. $\tau^*=0$ and $\sigma^*=0$ form a Nash equilibrium whenever $x\geq \log K$ with $V(x)=e^x$.
 In what follows, we assume




\bigskip

(A): $\psi(-1)<q.$

\bigskip

\noindent In that case, the Laplace exponent $\psi$ is well-defined on $[-1,\infty)$ and moreover  the L\'evy-Khintchine formula  can be extended to the interval $[-1,0)$ (see for instance Lemma 26.4 in  \cite{sa}). 
Without this condition the gain in the expectations in (\ref{val}) is infinity on the event $\{\tau=\sigma=\infty\}$.
\section{Main results.}
Below, in Theorems \ref{qsmall}--\ref{small-q} we give a qualitative and quantitative exposition
of the solution to (\ref{nash}). Before doing so, we need to give a brief reminder of
a class of special functions which appear commonly in connection with the
study of spectrally positive L\'evy processes. For each $p\ge 0$ we introduce the functions $W^{(p)}:\R\to [0,\infty)$ which are known to satisfy for all $x,y\ge 0$,
\begin{equation}\label{sfc}
\e\left[e^{-p\tau^-_{-x}}\ind_{\{\tau^-_{-x}<\tau^+_y\}}\right]=\frac{W^{(p)}(y)}{W^{(p)}(x+y)},
\end{equation}
where
\[
\tau^+_y:=\inf\{t>0:X_t>y\}\qquad\textrm{and}\qquad \tau^-_{-x}:=\inf\{t>0:X_t<-x\}
\]
(cf. Chapter 8 of Kyprianou \cite{ky}). In particular  $W^{(p)}(x)=0$ for all $x<0$ and further, it is known that on $(0,\infty)$ $W^{(p)}$ is almost everywhere differentiable, there is right continuity at zero and
\[
\int_0^{\infty}e^{-\beta x}W^{(p)}(x)\ud x=\frac{1}{\psi(\beta)-p}
\]
for all $\beta>\Phi(p)$, where $\Phi(p)$ is the largest root of the equation $\psi(\theta)=p$ (of which there are at most two).  For convenience, we write $W$ instead of $W^{(0)}$.

Associated to the functions $W^{(p)}$ are the functions $Z^{(p)}:\R\to [1,\infty)$ defined by
\[
Z^{(p)}(x)=1+p\int_{0}^x W^{(p)}(y)\ud y
\]
for $p\ge 0$. Together, the functions $W^{(p)}$ and $Z^{(p)}$ are collectively known as scale functions and predominantly appear in almost all fluctuation identities for spectrally positive L\'evy processes. For example, it is known that for all $y\in \R$
\begin{equation}\label{scf1}
\e\left[e^{-p\tau^+_{y}}\ind_{\{\tau^+_{y}<\infty\}}\right]=Z^{(p)}(y)-\frac{p}{\Phi(p)}W^{(p)}(y).
\end{equation}
We make the very mild assumption that $\Pi$ has no atoms when $X$ has paths of bounded variation. This suffices to deduce  (cf. \cite{savov}) that  $W^{(p)}\in C^1(0,\infty)$ and hence $Z^{(p)}\in C^2(0,\infty)$ and further, if $X$ has a Gaussian component they both belong to $C^2(0,\infty)$. It is also known that if $X$ has bounded variation with drift $d$, then $W^{(p)}(0+)=1/d$ and otherwise $W^{(p)}(0+)=0$. (Here and in the sequel we take the canonical representation of a bounded variation spectrally positive L\'evy processes $X_t=S_t-dt$ for $t\ge 0$ where $(S_t, t\ge 0)$ is a driftless subordinator and $d$ is a strictly positive constant which is referred to as the drift). Further, when $X$ has unbounded variation,
\begin{equation}\label{scfd}
W^{(p)\prime}(0+)=
2/b^2
\end{equation}
which is understood to be $+\infty$ when $b^2=0$.
Consider the exponential change of measure
\begin{equation}
\frac{\ud \p^{(\lambda)}}{\ud \p}\bigg|_{\mathcal{F}_t}=e^{-\lambda X_t-\psi(\lambda)t},\qquad \textrm{for } \lambda\ge -1.
\label{esscher}
\end{equation}
Under $\p^{(\lambda)}$, the process $X$ is still a spectrally positive L\'evy process and we mark its Laplace exponent and scale functions with the subscript ${\lambda}$. It holds that
\[
\psi_\lambda(\theta)=\psi(\lambda+\theta)-\psi(\lambda)
\]
for $\theta\ge 0$ and, by taking Laplace transforms, we find
\[
W^{(p)}_\lambda(x)=e^{-\lambda x}W^{(p+\psi(\lambda))}(x)
\]
for $p\ge 0$. The reader is otherwise referred to Chapter VII of Bertoin \cite{Be} or Chapter 8 of Kyprianou \cite{ky} for a general overview of  one-sided L\'evy processes and scale functions.

It turns out that the solution to the stochastic game can fall in four different regimes, depending on the value of the discount factor $q$. We remind the reader of the standing assumption (A).







\begin{theorem}\label{qsmall}
Suppose $q\leq \alpha/K$. Then a saddle point for the stochastic game   (\ref{nash}) is given by $\sigma^*=0,$ $\tau^*=\tau_{\log K}^+$. In particular $V(x)=K\vee e^x$ for all $x$.
\end{theorem}

\begin{theorem} \label{big-q} $\mbox{ }$
\begin{enumerate}
 \item[(i)]
As a function of $q$,
\[
a^*(q):=\frac{\alpha(\Phi(q)+1)}{\Phi(q)(q-\psi(-1)-\beta)}.
\]
is strictly monotone decreasing with $a^*(\beta+\psi(-1)+)=\infty$ and $a^*(\infty)=0$.
Define
\[
 q_0=\sup\Big\{q\in (0,\infty) :
a^*(q) < K
\Big\}
\]
 (noting that $q_0>\beta+\psi(-1)$ necessarily).
It holds that $q_0>\alpha/K$.

\item[(ii)]  For all $q\in [q_0, \infty)$
 a stochastic saddle point is given by the pair
\[
\tau^*=\inf\Big\{t\geq 0: X_t>\log a^*(q)\Big\}\quad\textrm{ and }\quad\sigma^*=\inf\Big\{t\geq 0: X_t>\log K\Big\}.
\]
In particular,
\[
V(x)=e^x +\frac{\alpha}{\Phi(q)}g(\log(a^*(q))-x),
\]
where
\[
g(z) = (\Phi(q)+1)\int_0^{z} e^{y-z} W^{(q)}(y){\rm d}y - \Phi(q)\int_0^{z}W^{(q)}(y){\rm d}y.
\]
\item[(iii)] For $q\in[q_0,\infty)$, there is smooth fit at $\log a^*(q)$ if and only if $X$ has paths of unbounded variation and otherwise there is continuous fit.
\end{enumerate}

\end{theorem}
\begin{theorem} \label{medium-q} \begin{enumerate}$\mbox{ }$
 \item[(i)] Assume that 
$b^2>0$.
 The set
\[
\left\{ q\in(0,q_0) : 0< K\frac{b^2}{2} +\frac{\alpha}{\Phi(q)}\left(\frac{K}{a^*(q)}-1\right)\right\}
\]
is an interval whose infimum we denote by $q_1$. In particular $q_1 > \alpha/K$.

\item[(ii)] When $b^2>0$ and  $q\in[q_1,q_0]$, a saddle point for  the stochastic game  (\ref{nash}) is given by
\[
\tau^*=\sigma^*=\inf\Big\{t\geq 0: X_t>\log K\Big\}.
\]
In particular,
\[
\begin{split}
V(x)&=e^x\left(1 + (q-\psi(-1))\int_0^{\log K-x} e^{y} W^{(q)}(y){\rm d} y -K\frac{q-\psi(-1)}{\Phi(q)+1}e^{-x}W^{(q)}(\log K-x)\right)\\
&\quad+\frac{W^{(q)}(\log K-x)}{K^{\Phi(q)}}\int_{-\infty}^{\log K} \big(\alpha+\beta e^{y}\big)e^{\Phi(q)y}\ud y-\int_x^{\log K} \big(\alpha+\beta e^{y}\big)W^{(q)}(y-x)\ud y.
\end{split}
\]
\item[(iii)] When $b^2>0$ and  $q\in[q_1, q_0]$ there is smooth fit at $\log K$ if and only if $q=q_0$ or $q=q_1$.
\item[(iv)] When $b^2 = 0$, then the strategies $\tau^*=\sigma^* = \inf\{t\geq 0 : X_t > \log K\}$ do not form  a stochastic saddle point when $q<q_0$. (In this case we define $q_1 = q_0$). 
\end{enumerate}
\end{theorem}

\begin{theorem}\label{small-q}
Suppose that $\alpha/K<q<q_1$.
\begin{enumerate}
 \item[(i)] The functional equation in $a$
\[
\begin{split}
\frac{q}{\Phi(q)}&=\frac{1}{K}\left(\frac{\alpha}{\Phi(q)}+\frac{\beta}{\Phi(q)+1}e^{a}\right)\\
&\qquad-\frac{1}{\Phi(q)}\int_{0}^\infty \Pi(\ud z+\log K-a)(1-e^{-\Phi(q)z})\\
&\qquad\quad+\frac{1}{\Phi(q)+1}\int_{0}^\infty \Pi(\ud z+\log K-a)e^z(1-e^{-(\Phi(q)+1)z})
\end{split}
\]
has a unique solution in $(-\infty, \log K)$ which we denote by $c^*(q)$.

\item[(ii)] A stochastic saddle point is given by the pair
\[
\tau^*=\inf\Big\{t\geq 0: X_t>\log K\Big\}\quad\textrm{ and }\quad\sigma^*=\inf\Big\{t\geq 0: X_t>c^*(q)\Big\}.
\]
In particular, for $x<c^*(q)$
\[
\begin{split}
V(x)&=K\Big(Z^{(q)}(c^*(q)-x)-\frac{q}{\Phi(q)}W^{(q)}(c^*(q)-x)\Big)\\
&+\int_{-\infty}^c \big(\alpha+\beta e^{y}\big)\Big(e^{-\Phi(q)(c^*(q)-y)}W^{(q)}(c^*(q)-x)-W^{(q)}(y-x)\Big)\ud y\\
&+e^{\Phi(q)(x-c^*(q))}\int_0^\infty e^{-\Phi(q)u}\Pi(\ud u) \int_{-\infty}^0 \ud y \Big(W^{(q)}(c^*(q)-x)-e^{-\Phi(q)y}W^{(q)}(c^*(q)-x+y)\Big)\\
&\hspace{5.3cm}\times e^{\Phi(q)(u+y)}\big(e^{c^*(q)+u+y}-K\big)\ind_{\{u+y+c^*(q)>\log K\}}
\end{split}
\]
and for $x\ge c^{*}(q)$
\[
V(x)=e^x\lor K.
\]
\item[(iii)] There is smooth fit at $c^*(q)$ if and only if $X$ has paths of unbounded variation and otherwise there is continuous fit.
\end{enumerate}

\end{theorem}
The order in which we present these statements  above (first $q\leq \alpha/K$, followed by $q\geq q_0$, $q\in[q_1,q_0]$ when $b^2>0$, and finally $q\in(\alpha/K,q_1)$) is convenient with regard to the dependency between their proofs. We also note that with the exception of Theorem \ref{medium-q}, the conclusions with regard to smooth versus continuous fit are consistent with existing results in the literature which generally exhibit smooth fit at boundary points of the stopping region if and only if that point is regular for the interior of the stopping region (see for example \cite{aky} and \cite{PS}). In our case, thanks to stationary and independent increments of $X$, this boils down to the point $0$ being regular for $(0,\infty)$, which also corresponds to the case that $X$ has unbounded variation for the special case of spectrally positive L\'evy processes.

The remainder of this paper is dedicated to proving these theorems and is structured as follows. In the next section we state and prove a Lemma which will be repeatedly used to implement proofs on the basis of `guess and verify' such as is common with solving optimal stopping problems. Thereafter we prove the four main theorems above in the order that they are stated.

\section{Guess and verify}
Following classical ideas in optimal stopping, we verify that a candidate solution solves the stochastic game (\ref{nash}) by checking certain associated bounds and martingales properties. Specifically, we use the following verification lemma which is of a similar form to Lemma 5 in \cite{bky}.
\begin{lemma}[Verification Lemma]  Fix $x\in\mathbb{R}$. Suppose that  $\tau^*\in \mathcal{T}_{0,\infty}$ and  $\sigma^*\in \mathcal{T}_{0,\infty}$ are candidate optimal strategies for the stochastic game (\ref{stgame}) such that
\begin{equation}\label{condvl}
\sup_{\sigma\in \mathcal{T}_{0,\infty}}e^{-q\sigma +X_{\sigma}}\ind_{\{\sigma\le \tau^*\}}
\end{equation}
has finite mean under $\mathbb{P}_x$. Let
\[
V^*(x)=\ex\Big(L_{\tau^*}\ind_{\{\tau^*<\sigma^*\}}+U_{\sigma^*}\ind_{\{\sigma^*\leq \tau^*\}}\Big).
\]
Then  $(\tau^*, \sigma^*)$ is a stochastic saddle point  of (\ref{nash}) with value $V^*$ if
\begin{itemize}
\item[(i)] $V^*(x)\geq e^x$,
\item[(ii)] $V^*(x)\leq e^x\lor K$,
\item[(iii)] $V^*(X_{\tau^*})=e^{X_{\tau^*}}$ almost surely on $\{\tau^*<\infty\}$,
\item[(iv)] $V^*(X_{\sigma^*})=e^{X_{\sigma^*}}\lor K$ almost surely on $\{\sigma^*<\infty\}$,
\item[(v)] the process
\[
\left(e^{-q(t\land \tau^*)}V(X_{t\land \tau^*})+\int_0^{t\land \tau^*} e^{-qs}(\alpha+\beta e^{X_s})\ud s,\, t\ge 0\right)
\]
is a right continuous submartingale and
\item[(vi)] the process
\[
\left(e^{-q(t\land \sigma^*)}V(X_{t\land \sigma^*})+\int_0^{t\land \sigma^*} e^{-qs}(\alpha+\beta e^{X_s})\ud s,\, t\ge 0\right)
\]
is a right continuous supermartingale.
\end{itemize}
\end{lemma}
\noindent{\it Proof:} Define
\[
H^q_{\tau,\sigma}:=L_{\tau}\ind_{\{\tau< \sigma\}}+U_{\sigma}\ind_{\{\sigma\le\tau\}},
\]
where $\tau$ and $\sigma$ are stopping times.  Since we have assumed that $q>\psi(-1)\vee 0$, we have that
\[
H^q_{\infty,\infty}=\int_0^\infty e^{-qs}(\alpha+\beta e^{X_s})\ud s.
\]
From the supermartingale property $(vi)$, Doob's optional sampling theorem, $(iv)$ and $(i)$ we know that for any stopping time $\tau$ and $t\ge 0$,
\[
\begin{split}
V^*(x)&\ge \ex\left(e^{-q(t\land \tau\land \sigma^*)}V^*(X_{t\land \tau\land \sigma^*})+\int_0^{t\land \tau\land \sigma^*}e^{-qs}(\alpha+\beta e^{X_s})\ud s\right)\\
&\ge \ex\bigg(e^{-q(t\land \tau)+X_{\tau\land t}}\ind_{\{t\land \tau<\sigma^*\}}+e^{-q\sigma^*}(e^{X_{\sigma^*}}\lor K)\ind_{\{\sigma^*\le t\land \tau\}}\\
&\qquad\quad+\int_0^{t\land \tau\land \sigma^*}e^{-qs}(\alpha+\beta e^{X_s})\ud s\bigg).
\end{split}
\]
It follows from Fatou's lemma by taking $t$  to $\infty$, that
\[
V^*(x)\ge \ex(H^q_{\tau, \sigma^*}).
\]
Now using $(v)$, Doob's optional sampling theorem, $(iii)$ and $(ii)$ we have  for any stopping time $\sigma$ and $t\ge 0$,
\[
\begin{split}
V^*(x)&\le \ex\bigg(e^{-q\tau^*+X_{\tau^*}}\ind_{\{ \tau^*<\sigma\land t\}}\bigg)+\ex\bigg(\int_0^{\tau^*\land \sigma\land t}e^{-qs}(\alpha+\beta e^{X_s})\ud s\bigg)\\
&\qquad\quad+\ex\bigg(e^{-q(\sigma\land t)}(e^{X_{\sigma\land t}}\lor K)\ind_{\{\sigma\land t\le \tau^*\}}\bigg).
\end{split}
\]
Taking limits as $t$ goes to $\infty$ and applying the monotone convergence theorem for the first two  terms on the right-hand side and the dominated convergence theorem for the last term on the right-hand side (see (\ref{condvl})), we have
 \[
V^*(x)\le \ex(H^q_{\tau^*, \sigma}),
\]
and hence $(\tau^*, \sigma^*)$ is a saddle point to (\ref{nash}). \hfill$\square$

\section{Proof of Theorem \ref{qsmall}}

Suppose $q\leq \alpha/K$. We claim that the process $(Z_t, \,t \geq 0)$ defined by
\[Z_t=\left(e^{-q(t\wedge\tau_{\log K}^+)}(K\vee e^{X_{t\wedge\tau_{\log K}^+}})+\int_0^{t\wedge\tau_{\log K}^+} e^{-qs}(\alpha+\beta e^{X_s})\ud s,\, t\ge 0\right) \]
is a submartingale.
Indeed, when  $x<\log K$, we have on 
$\{t<\tau^+_{\log K}\}$
\[
dZ_t = [(\alpha-qK)e^{-qt} + \beta e^{-qt+ X_t}]{\rm d}t
\]
and 
\[
Z_{\tau^+_{\log K}} - Z_{\tau^+_{\log K} - } = e^{-q \tau^+_{\log K}}(e^{X_{\tau^+_{\log K}}} - K),
\]
showing that, as $\beta>0$, $Z$ is an adapted, strictly increasing process, i.e. a submartingale.


We may now invoke the Verification Lemma, since the other properties are automatically satisfied by taking $\sigma^*=0$. Note in particular that the condition (\ref{condvl}) is automatically satisfied since
\[
\sup_{\sigma\in \mathcal{T}_{0,\infty}}e^{-q\sigma + X_{\sigma}}\ind_{\{\sigma\le \tau^+_{\log K}\}}
\leq e^{-q\tau^+_{\log K} + X_{\tau^+_{\log K}}}\ind_{\{\tau^+_{\log K}<\infty\}}+K
\]
and
\[
\mathbb{E}_x(e^{-q \tau^+_{\log K } +X_{\tau^+_{\log K}}}\ind_{\{\tau^+_{\log K}<\infty\}} ) =
e^x\mathbb{E}^{(-1)}_x(e^{-(q-\psi(-1) )\tau^+_{\log K } } ) <e^x.
\]


\section{Proof of Theorem \ref{big-q}}\label{proof-big-q}

The basis of the proof of Theorem \ref{big-q} is the assumption that the optimal strategies take the form
$\sigma^* = \inf\{ t>0 : X_t > \log K\}$ and $\tau^* = \inf\{t> 0 : X_t > y^*\}$ for some optimally chosen $y^*$. On this basis, establishing the value function in the Stackelberg equilibrium, $V$, would boil down to computing $H_{y^*}$ where for any $-\infty<y \leq \log K$,
$H_y(x):=\mathbb{E}_x\big(L_{\tau^+_y}\big)$, that is to say,
\begin{equation}\label{ident1}
H_y(x) :=\mathbb{E}_x\Big(e^{-q\tau^+_y+X_{\tau^+_y}}\Big)+\int_0^{\infty}\mathbb{E}_x \Big(e^{-qs}(\alpha+\beta e^{X_s})\ind_{\{s\le \tau^+_y\}}\Big)\ud s.
\end{equation}
We thus proceed by evaluating the above expression in terms of scale functions, then we choose the value of $y^*$ by blindly applying the principle of smooth and continuous fit respectively to the cases that $X$ has paths of unbounded and bounded variation and finally we verify that the established strategy is indeed optimal with the help of the Verification Lemma.

With the help of the exponential change of measure, (\ref{esscher}),  and (\ref{scf1}), the first term of the right-hand side of the above expression for $H_y$ satisfies
\[
\begin{split}
\mathbb{E}_x\Big(e^{-q\tau^+_y+X_{\tau^+_y}}\Big)
&=e^x\mathbb{E}\Big(e^{-q\tau^+_{y-x}+X_{\tau^+_{y-x}}}\Big)\\
&=e^x\mathbb{E}^{(-1)}\Big(e^{-(q-\psi(-1))\tau^+_{y-x}}\Big),\\
&=e^x\left(Z_{-1}^{(q-\psi(-1))}(y-x)-\frac{q-\psi(-1)}{\Phi_{-1}(q-\psi(-1))}e^{y-x}W^{(q)}(y-x)\right)\\
&=e^x\left(1 +(q-\psi(-1))\int_0^{y-x}e^{y}W^{(q)}(z){\rm d}z-\frac{q-\psi(-1)}{\Phi_{-1}(q-\psi(-1))}e^{y-x}W^{(q)}(y-x)\right).
\end{split}
\]
On the other hand from Theorem 8.7 in \cite{ky}, the second term of the right-hand side of (\ref{ident1}) satisfies
\[
\begin{split}
\int_0^{\infty}\mathbb{E}_x \Big(e^{-qs}(\alpha+\beta e^{X_s})\ind_{\{s\le \tau^+_y\}}\Big)\ud s&=\int_0^{\infty}\widehat{\mathbb{E}} \Big(e^{-qs}(\alpha+\beta e^{x-X_s})\ind_{\{s\le \tau^-_{x-y}\}}\Big)\ud s\\
&=\int_0^{\infty}\big(\alpha+\beta e^{y-z}\big)\Big(e^{-\Phi(q)z}W^{(q)}(y-x)-W^{(q)}(y-x-z)\Big)\ud z\\
&=\int_{-\infty}^y \big(\alpha+\beta e^{z}\big)\Big(e^{-\Phi(q)(y-z)}W^{(q)}(y-x)-W^{(q)}(z-x)\Big)\ud z,
\end{split}
\]
where $\widehat{\p}$ denotes the law of the dual process $\widehat{X}=-X$.
Finally noting that $\Phi_{-1}(q-\psi(-1))=\Phi(q)+1$, we get
\[
\begin{split}
H_y(x)&=e^x\left(1 +(q-\psi(-1))\int_0^{y-x}e^{z}W^{(q)}(z){\rm d}z-\frac{q-\psi(-1)}{\Phi(q)+1}e^{y-x}W^{(q)}(y-x)\right)\\
&\quad+\int_{-\infty}^y \big(\alpha+\beta e^{z}\big)e^{-\Phi(q)(y-z)}W^{(q)}(y-x)\ud z-\int_x^y \big(\alpha+\beta e^{z}\big)W^{(q)}(z-x)\ud z.
\end{split}
\]
We also see in particular, making use of the fact that $W^{(q)}(0-)=0$ and $Z^{(q)}(0)=1$, that
\[
H_y(x)=e^x
\]
for all $x>y$.

\bigskip

Having expressed $H_y$ in terms of scale functions, we now turn our attention to making the choice of $y^*$ using the principle of smooth and continuous fit.

\bigskip

{\bf Bounded variation and continuous fit:} In this case it is known that $W^{(q)}(0+)>1/d$  where $d>0$ is the drift term of the process $X$. It follows that

\begin{eqnarray}
H_{\log{a}}(\log{a}-)&=&a\left(1-\frac{q-\psi(-1)}{\Phi(q)+1}W^{(q)}(0+)\right)+W^{(q)}(0+)a^{-\Phi(q)}\int_{-\infty}^{\log a} \big(\alpha+\beta e^{y}\big)e^{\Phi(q)y}\ud y \notag\\
&=&a+aW^{(q)}(0+)\left(\frac{\alpha}{ a \Phi(q)}+\frac{\beta}{\Phi(q)+1} -\frac{q-\psi(-1)}{\Phi(q)+1}W^{(q)}(0+)\right).
\label{fit}
\end{eqnarray}

In order to avoid a discontinuity at $a$ we choose it equal to the value $a^*$ which satisfies
\[
\frac{q-\psi(-1)}{\Phi(q)+1}=\frac{\alpha}{a\Phi(q)}+\frac{\beta}{\Phi(q)+1}.
\]
Note that this is equivalent to requiring that
\begin{equation}
a^*= \frac{\alpha(\Phi(q)+1)}{\Phi(q)(q-\psi(-1)-\beta)}
\label{a*}
\end{equation}
providing $q>\psi(-1)+\beta$. In order to respect the requirement that $a^*\leq  K$ we also need to check how the function $a^*= a^*(q)$ varies with $q$. To this end, note that
\[
\frac{d}{dq}a^*(q)=\frac{-(q-\psi(-1)-\beta)\alpha\Phi'(q)-\alpha(\Phi(q)+1)\Phi(q)}{(\Phi(q)(q-\psi(-1)-\beta))^2} <0,
\]
hence $a^*(\cdot)$ is strictly decreasing.  Note also that
\[
\lim_{q\downarrow \beta+\psi(-1)}a^*(q)=\infty\qquad\textrm{ and }\qquad\lim_{q\to \infty}a^*(q)=0,
\]
which implies the existence of a unique $q_0>\beta+\psi(-1)$ such that $a^*(q_0)=K$.  Note that it also turns out that $q_0> \alpha/K$ on account of the fact that for $q\leq \alpha/K$
\[
a^*(q) \geq K \frac{q}{\Phi(q)}\frac{\Phi(q)+1}{q-\psi(-1)- \beta}\geq K \frac{q}{\Phi(q)}\frac{\Phi(q)+1}{q-\psi(-1)} = K \mathbb{E}(e^{\overline{X}_{\mathbf{e}_q}})> K
\]
where in the equality we have appealed to the well-known identity for one of the Wiener--Hopf factors  of $X$ (cf. Chapter 8 of \cite{ky}).

\bigskip

{\bf Unbounded variation and smooth fit:} In this case it is known that $W^{(q)}(0+)=0$ and hence in the above analysis one sees that $H_{\log a}(\log a-)=a = H_{\log a}(\log a+)$. In that case, the principle of smooth fit can be implemented and we insist on there being no discontinuity in $H_{\log a}'$ at $\log a$.
We have
\begin{eqnarray}
H_{\log a}'(x)&=& e^xZ^{(q-\psi(-1))}_{-1}(\log(a)-x) - e^x (q-\psi(-1)) W^{(q-\psi(-1))}_{-1}(\log(a)-x) \notag\\
&&+ a\frac{q-\psi(-1)}{\Phi(q)+1}W^{(q)\prime}(\log(a)-x) -W^{(q)\prime}(\log(a)-x) a^{-\Phi(q)}\int_{-\infty}^{\log a} \big(\alpha+\beta e^{y}\big)e^{\Phi(q)y}\ud y \notag\\
&&+(\alpha +\beta e^x)W^{(q)}(0+). \label{prime}
\end{eqnarray}
Recall that $W^{(q)}(0+)=0$ and that $W^{(q)\prime}(0+)=2/b^2$ which should be interpreted as $+\infty$ in the case that the Gaussian coefficient $b^2 =0$. We find
\begin{equation}
H_{\log a}'(\log a-)=a +\left[ a\frac{q-\psi(-1)}{\Phi(q)+1}-a^{-\Phi(q)}\int_{-\infty}^{\log a} \big(\alpha+\beta e^{y}\big)e^{\Phi(q)y}\ud y\right]W^{(q)\prime}(0+).
\label{uselater}
\end{equation}
In order to obtain the smooth fit $H_{\log a}'(\log a+)=a$ we must thus have that
\[
a\frac{q-\psi(-1)}{\Phi(q)+1}=a^{-\Phi(q)}\int_{-\infty}^{\log a} \big(\alpha+\beta e^{y}\big)e^{\Phi(q)y}\ud y,
\]
which, after a simple integration on the right hand side,
gives the same expression of $a^*$ as in the bounded variation case. The same bounds on $q_0$ are thus still applicable in this case too.

\bigskip

In both cases, we obtain our candidate value function
\begin{eqnarray*}
H_{\log a^*}(x)&=&e^x Z^{(q-\psi(-1))}_{(-1)}(\log(a^*)-x)-\int_x^{\log a^*} \big(\alpha+\beta e^{y}\big)W^{(q)}(y-x)\ud y\\
&=&e^x\left(1 + (q-\psi(-1)) \int_0^{\log(a^*)-x} e^yW^{(q)}(y)dy\right)-\int_0^{\log(a^*)-x} \big(\alpha+\beta e^{y+x}\big)W^{(q)}(y)\ud y.
\end{eqnarray*}

\bigskip

We now proceed to verify our candidate solution when $q\geq q_0$.  That is to say, we shall verify that
\[
\tau^*=\tau^+_{\log a^*}, \qquad \sigma^*=\tau^+_{\log K} \mbox{ and } V^*(x) = H_{\log a^*}(x)
\]
fulfill the conditions of the Verification Lemma. Note in particular that $\tau^*\leq \sigma^*$.


\bigskip

{\bf Submartingale and supermartingale  properties:}  To this end, note  that  (\ref{ident1}) together with an application of the Markov property gives us for all $t\geq 0$,
\begin{eqnarray*}
\Lambda_t &: =& \mathbb{E}_x\left[\left.
e^{-q\tau^+_{\log a^*}}V^*(X_{ \tau^+_{\log a^*}})
+\int_0^{\tau^+_{\log a^*}} e^{-qs}(\alpha+\beta e^{X_s})\ud s \right|\mathcal{F}_t\right]  \\
&&=e^{-q(t\wedge\tau^+_{\log a^*})}V^*(X_{t\wedge \tau^+_{\log a^*}})
+\int_0^{t\wedge \tau^+_{\log a^*}} e^{-qs}(\alpha+\beta e^{X_s}) \ud s.
\label{useSMP}
\end{eqnarray*}
That is to say, $\Lambda = (\Lambda_t : t\geq 0)$  is a martingale.
This confirms the submartingale property (v) in the Verification Lemma.

An easy computation shows that 
\[
 V^{*\prime\prime}(x)= e^xZ^{(q-\psi(-1))}_{-1}(\log (a^*)-x)  - a^*(q-\psi(-1))[W^{(q)}(\log(a^*)-x) + W^{(q)\prime}(\log(a^*)-x)]
\]
and hence $V^*$ belongs to $C^2(-\infty, \log a^*)$. Moreover, the latter conclusion is sufficient to show that $\Gamma V^*(x)$ is continuous on $(-\infty, \log a^*)$ where $\Gamma$ is the infinitesimal generator of $X$, and in particular,
\[
 \Gamma V^*(x) = \mu V^{*\prime}(x) + \frac{b^2}{2}V^{*\prime\prime}(x)
+ \int_{(0,\infty)} (V^*(x+y) - V^*(x) - yV^{*\prime}(x)\ind_{\{y<1\}})\Pi({\rm d}y).
\]
 (See for example the argument in Lemma 4.1 of \cite{krs}). For any $-n\leq x\leq a<\log a^*$, where $n\in \mathbb{N},$ the aforementioned facts concerning smoothness and continuity allow us to apply
 It\^o's formula 
to $\Lambda$, but stopped at  $\tau^-_{-n}$, where $\tau^-_{-n} = \inf\{t>0: X_t <-n\}$, and deduce that
\begin{equation}
\Lambda_{t \wedge \tau^-_{-n}} = V^*(x) + \int_0^{t\wedge \tau^-_{-n}\wedge\tau^+_a} e^{-qs}[
(\Gamma - q)V^*(X_s)  + (\alpha + \beta e^{X_s}) ]{\rm d}s + m_t
\label{nodrift}
\end{equation}
where
\begin{eqnarray*}
m_t &=&\frac{b^2}{2}\int_0^{t \wedge \tau^-_{-n}\wedge \tau^+_a}e^{- qs}V^{*\prime}(X_s){\rm d}B_s
+  \int_0^{t \wedge \tau^-_{-n}\wedge \tau^+_a}e^{- qs}V^{*\prime}(X_s){\rm d}X^{(1)}_s\\
&&+ \sum_{s\leq t \wedge \tau^-_{-n}\wedge \tau^+_a} e^{-qs} \left[
V^*(X_s)- V^*(X_{s-}) - \Delta X_s V^{*\prime}(X_{s-})\ind_{\{\Delta X_s <1\}}
\right]\\
&&- \int_0^{t \wedge \tau^-_{-n}\wedge \tau^+_a} e^{-qs}  \int_{(0,\infty)} (V^*(X_{s-}+y) - V^*(X_{s-}) - yV^{*\prime}(X_{s-})\ind_{\{y<1\}})\Pi({\rm d}y) {\rm d}s
\end{eqnarray*}
is a local martingale such that $B$ is the Gaussian component in $X$ and $X^{(1)}$ is the martingale
part of X consisting of compensated jumps of size strictly less than unity. In fact, thanks to the boundedness of $V^{*\prime}$ and $\Gamma V^{*}$ on $[-n,a]$, the process $\{m_t: t\geq 0\}$ is a martingale.
The latter, together with the fact that $\Lambda$ is a martingale, implies that the drift term in (\ref{nodrift})
must
almost surely be equal to zero.
 Taking expectations and writing $R^{(q)}(x, {\rm d}y ; a, -n)$ for the $q$-resolvent measure of  the process $X$ when issued from $x$ and killed on first entry into $(-\infty, -n)\cup (a,\infty)$ we have for all $-n\leq x\leq a<\log a^*$,
\[
 \int_{[-n,a]} [(\Gamma - q)V^*(y)  + (\alpha + \beta e^{y}) ] R^{(q)}(x, {\rm d}y ; a, -n) = 0.
\]
As $R^{(q)}(x, {\rm d}y ; a, -n)$ is absolutely continuous with respect to Lebesgue measure with a strictly positive density in $(-\infty, 0)$ (cf. Chapter 8 of Kyprianou \cite{ky}) it follows that
\begin{equation}
(\Gamma - q)V^*(x)  + (\alpha + \beta e^{x})  = 0
\label{iscts}
\end{equation}
for Lebesgue almost every $x< \log a^*$. The latter can be upgraded to every $x< \log a^*$ as the left hand side of (\ref{iscts}) is continuous.
It is also trivial to check that  $V^*(x) = e^x$ on $(\log a^*,\infty)$ and hence it follows from $q>\psi(-1)+\beta$ that
\[
(\Gamma - q)V^*(x)  + (\alpha + \beta e^x)= (\psi(-1)-q+\beta)e^x + \alpha
\leq (\psi(-1)-q+\beta)a^* + \alpha\leq 0
\]
on $(\log a^*,\infty).$

Next, note that it is straightforward to 
 see that  $V^*$ is twice continuously differentiable on $(-\infty, \log a^*)\cup(\log a^*,\infty)$ with the existence of a left and right derivative at $\log a^*$. 
 We may thus apply the Meyer--It\^o formula (cf. Theorem 70 of Protter \cite{protter}) to the process $V(X_{t\wedge \tau^+_{\log K}})$ and then integrate by parts to obtain, in a similar vein to (\ref{nodrift}), that
\begin{eqnarray*}
\lefteqn{e^{-q(t\wedge \tau^+_{\log K})} V^*(X_{t\wedge \tau^+_{\log K}}) + \int_0^{t\wedge \tau^+_{\log K}} e^{-qs}(\alpha + \beta e^{X_s}) {\rm d}s}\\
 &&= V^*(x)
+ \int_0^{t\wedge \tau^+_{\log K}} e^{-qs}[
(\Gamma - q)V^*(X_s)  + (\alpha + \beta e^{X_s}) ]{\rm d}s \\
&&+\frac{1}{2}\int_0^{t\wedge \tau^+_{\log K}} e^{-qs}(V^{*\prime}(\log a^* +) - V^{*\prime}(\log a^* -) ) {\rm d}\ell_s + M_t
\end{eqnarray*}
where $M: = (M_t: t\geq 0)$ is a local martingale and $\ell: = (\ell_t: t\geq 0)$ is the semi-martingale local of $X$ at $\log a^*$. Note that when $b^2=0$, the final integral is identically zero owing to the fact that the local time process $\ell$ is also identically zero and otherwise, when $b^2>0$, the final integral is still identically zero thanks to smooth pasting. Note also that although the quantity $(\Gamma - q)V^*(x)  + (\alpha + \beta e^{x}) $ is not defined at $x=\log a^*$, this is not a problem in the context of the above calculus as the Lebesgue measure of the time that the process $X$ spends at $\log a^*$ is zero. 

Recalling  that $(\Gamma - q)V^*(x)  + (\alpha + \beta e^{x})\leq 0 $ on $(-\infty, \log a^*)\cup(\log a^*,\infty)$,  by taking expectations with the help of a suitable localizing sequence of stopping times $\{T_n : n\geq 1\}$ for $M$, Fatou's lemma and monotone convergence, we obtain
\begin{eqnarray*}
\lefteqn{ \mathbb{E}_x\left[e^{-q(t\wedge \tau^+_{\log K}) }     V^*(X_{t\wedge\tau^+_{\log K}})
+\int_0^{t\wedge \tau^+_{\log K}}e^{-qs}(\alpha+\beta e^{X_s}) \ud s
\right]} \\
&&\leq \lim_{n\uparrow\infty}  \mathbb{E}_x\left[e^{-q(t\wedge T_n \wedge \tau^+_{\log K}) }     V^*(X_{t\wedge T_n \wedge\tau^+_{\log K}})
+\int_0^{t\wedge T_n \wedge  \tau^+_{\log K}}e^{-qs}(\alpha+\beta e^{X_s}) \ud s
\right]\\
&&\leq  V^*(x)
+  \lim_{n\uparrow\infty} \mathbb{E}_x\left[  \int_0^{t \wedge T_n \wedge \tau^+_{\log K}} e^{-qs}[
(\Gamma - q)V^*(X_s)  + (\alpha + \beta e^{X_s}) ]{\rm d}s\right]\\
&&\leq V^*(x) .
\end{eqnarray*}
The last inequality above together with the Markov property is sufficient to deduce the supermartingale property (vi) in the Verification Lemma. Note that right continuity follows immediately from the continuity of $V^*$ and the fact that $X$ has c\`adl\`ag paths.

\bigskip

{\bf Lower and Upper bounds:} The bounds (i) and (ii) in the Verification Lemma can be deduced directly from the expression for $V^*$. To this end, write
\begin{eqnarray*}
V^*(x) 
&=&e^x + \frac{\alpha}{\Phi(q)} g(\log(a^*) - x),
\end{eqnarray*}
where
\[
 g(z) = (\Phi(q)+1)\int_0^{z} e^{y-z} W^{(q)}(y){\rm d}y - \Phi(q)\int_0^{z}W^{(q)}(y){\rm d}y.
\]
Note that $g(0)=0$. Since $V^*(x) = e^x$ for all $x\geq \log a^*$, we have the required lower bound for $V^*$ if we can prove that $g'(z)>0$ for all $z>0$.
To this end we differentiate and find that
\begin{eqnarray*}
 g'(z) &=& e^{-z}\left[e^{z}W^{(q)}(z) - (\Phi(q)+1)\int^z_0 e^y W^{(q)}(y){\rm d}y\right]\\
&=& e^{-z}\left[ \mathcal{W}^{(p)}(z) - \varphi(p)\int_0^z \mathcal{W}^{(p)}(y){\rm d}y \right],
\end{eqnarray*}
where $p=q-\psi(-1)$, $\mathcal{W}^{(p)}(z) = e^{z}W^{(q)}(z) = W_{-1}^{q-\psi(-1)}(z)$ and
 \begin{eqnarray*}
  \varphi(p)&=& \sup\{\theta\geq 0 : \psi_{-1}(\theta) = p\}\\
&=& \sup\{\theta\geq 0 : \psi(\theta-1) - \psi(-1) = q- \psi(-1)\}\\
&=&\sup\{\theta\geq 0 : \psi(\theta-1)  = q\}\\
&=& \Phi(q)+1.
 \end{eqnarray*}
Finally, to show that $g'(z)>0$
we note from (8.20) of Kyprianou (2006) that
\[
 0<\frac{\varphi(p)}{p}\mathbb{P}^{(-1)}(-\underline{X}_{\mathbf{e}_p}\leq z) = \mathcal{W}^{(p)}(z) - \varphi(p)\int_0^z \mathcal{W}^{(p)}(y){\rm d}y,
\]
where $\mathbf{e}_p$ is an exponentially distributed random variable which is independent of $X$ and has parameter $p$.

For the upper bound on $V^*$ it suffices to show in a similar vein to the lower bound that $V^{*\prime}(x)\geq 0$. Calculations in the spirit of the ones above show that
\begin{eqnarray*}
V^{*\prime}(x) &=& \frac{e^{x}}{a^*}\left[a^* - \frac{\alpha}{\Phi(q)} \frac{\Phi(q)+1}{q-\psi(-1)} \mathbb{P}^{(-1)}(-\underline{X}_{\mathbf{e}_p}\leq x-\log a^*)\right]\\
&\geq &\frac{e^x}{a^*}\frac{\alpha}{\Phi(q)} \frac{\Phi(q)+1}{q-\psi(-1)}[1 - \mathbb{P}^{(-1)}(-\underline{X}_{\mathbf{e}_p}\leq x-\log a^*)]\\
&\geq &0,
\end{eqnarray*}
where we have made use of (\ref{a*}).

\bigskip

{\bf Stopped values:} Note that since $V^*(x) = e^x$ for $x\geq \log a^*$ both conditions (iii) and (iv) are automatically satisfied.

\bigskip

Having now checked properties (i)-(vi) of the Verification Lemma, and noting that the justification for (\ref{condvl}) is the same as in the proof of Theorem \ref{qsmall}, we may conclude that the proposed triple $(\tau^*, \sigma^*, V^*)$ is a stochastic saddle point. \hfill$\square$

 \section{Proof of Theorem \ref{medium-q}}\label{middlecase}
The proof of Theorem \ref{medium-q} relies on the following optimal stopping problem.
 Recall that for $q>0$
\[
U_t=e^{-qt}(e^{X_t}\lor K)+\int_0^t e^{-qs}(\alpha+\beta e^{X_s})\ud s.
\]
\begin{lemma} \label{lemt3}Let $\alpha/ K<q\leq q_0$.
Define the function $w(x):\R\to \R^+$ by
\begin{equation}\label{opstop}
w(x):=\inf_{\sigma\in \mathcal{T}_{0,\infty}}\ex\Big[U_{\sigma\land \tau^+_{\log K}}\Big].
\end{equation}
Then $w$ has the following properties,
\begin{enumerate}
\item[(i)] $w$ is  non-decreasing,
\item[(ii)] $w(x)\leq e^x\lor K$ for $x\in \R$,
\item[(iii)] there exists a $c^*\leq \log K$ such that
\[
w(x)=\ex\left[e^{-q\tau^+_{c^*}}\Big(e^{X_{ \tau^+_{c^*}}}\lor K\Big)+\int_0^{ \tau^+_{c^*}} e^{-qs}(\alpha+\beta e^{X_s})\ud s\right],
\]
\item[(iv)]  w is continuous in $x$ and $w(c^*)=K$,
\item[(v)] $w(x)\ge e^x$ for $x\in \R$,
\item[(vi)] $w(X_{\tau^+_{\log K}})=e^{X_{\tau^+_{\log K}}}$ almost surely on $\{\tau^+_{\log K}<\infty\}$
\item[(vii)] $w(X_{\tau^+_{c^*}})=e^{X_{\tau^+_{c^*}}}\lor K$ almost surely on $\{\tau^+_{c^*}<\infty\}$
\item[(viii)] the process
\[
\left(e^{-q(t\land \tau^+_{\log K})}w(X_{t\land \tau^+_{\log K}})+\int_0^{t\land \tau^+_{\log K}} e^{-qs}(\alpha+\beta e^{X_s})\ud s,\, t\ge 0\right),
\]
is a right continuous submartingale and
\item[(ix)] the process
\[
\left(e^{-q(t\land \tau^+_{c^*})}w(X_{t\land\tau^+_{c^*}})+\int_0^{t\land \tau^+_{c^*}} e^{-qs}(\alpha+\beta e^{X_s})\ud s,\, t\ge 0\right),
\]
is a right continuous supermartingale.
\end{enumerate}
\end{lemma}
\noindent{\it Proof:} $(i)$ Denote $X^*_t=X_{t\land \tau^+_{\log K}}$ for all $t\ge 0$, and introduce the additive functional
\[
A_t:=\int_0^{t\land \tau^+_{\log K}} e^{-qs}(\alpha+\beta e^{X_s})\ud s,\qquad \textrm{ for all } \quad t\ge 0.
\]
Then the process $Z:=(Z_t,t\ge0)$  given by
\[
Z_t:=(t,A_t,X^*_t)\qquad \textrm{ for all } \quad t\ge 0,
\]
is Markovian and starts from  $(0,0,x)$ under the measure $\p_x$.  For $(t, i, x)\in \R^2_+\times\R$ denote by $\mathrm{P}_{(t,i,x)}$ the law of $Z$ when it started at $(t, i, x)$.  Thus the optimal stopping problem $(\ref{opstop}$)
 reads as follows
 \[
 w(x):=W(0,0,x)=\inf_{\sigma\in \mathcal{T}_{0,\infty}}\mathrm{E}_{(0,0,x)}\Big[F(\sigma\land\tau^+_{\log K}, A_{\sigma},X^*_\sigma)\Big],
 \]
 where $F(t,i,x)=e^{-qt}(e^x\lor K)+i$. Since $F:\R^2_+\times\R\to \R_+$ is continuous and $X^*$ is quasi-left continuous we can deduce that $w$ is upper semicontinuous. Furthermore, we have
 \[
 \begin{split}
 \mathrm{E}_{(0,0,x)}\left[\sup_{t\ge 0}F(t\land\tau^+_{\log K}, A_t,X^*_t)\right]&\leq \ex\left[\int_0^{\infty} e^{-qs}(\alpha+\beta e^{X_s})\ud s\right]+\ex\Big[e^{-q\tau^+_{\log K}+X_{\tau^+_{\log K}}}\Big]+K\\
&\le K+\frac{\alpha}{q}+\frac{\beta}{q-\psi(-1)}+e^x Z^{(q-\psi(-1))}_{-1}(\log K-x)<\infty,
\end{split}
 \]
 so we can apply a variant of Theorem 3.3 on p.127 of Shiryaev \cite{Shir} (see also Corollary 2.9 on p.46 of Peskir and Shiryaev \cite{PS}) to conclude that
 \[
 \tau_D=\inf\{t\ge 0:Z_t\in D\},
 \]
 where $D=\big\{(t,i,x)\in\R^2_+\times\R : W(t,i,x)=F(t,i,x)\big\}$, is
  an optimal stopping time. Note that for all $(t,i,x)\in \R^2_+\times\R$, the following identity holds
  \[
  W(t,i,x)=e^{-qt}W(0,0,x)+i,
  \]
and thus we deduce that $D=\{x\in \R : w(x)=e^x\lor K\}$ and  $ \tau_D=\tau^+_{\log K}\land \inf\{t\ge 0: X_t\in D\}$.

In what follows, if $\varsigma$ is a stopping time for $X$ we shall write $\varsigma(x)$ to show the dependence of the stopping time on the value of $X_0=x$. Similarly, we denote
\[
U^{(x)}_t=e^{-qt}(e^{X_t+x}\lor K) +\int_0^{t} e^{-qs}(\alpha+\beta e^{X_s +x})\ud s,\qquad t\geq 0.
\]
For $y\ge x$, we have that   $U^{(y)}_t\ge U^{(x)}_t$ for all $t\ge 0$  and thus, also appealing to the definition of $w$ as an infimum,
\[
\begin{split}
w(x)-w(y)=\e\Big[U^{(x)}_{\tau_D(x)}-U^{(y)}_{\tau_D(y)}\Big]\leq \e\Big[U^{(x)}_{\tau_D(y)}-U^{(y)}_{\tau_D(y)}\Big]\leq 0,
\end{split}
\]
which implies that $w$ non-decreasing.\\
$(ii)$ This property follows directly from the definition of $w$ as an infimum and taking for instance  the stopping time  $\sigma=0$.\\
$(iii)$ Recall that $w$ is upper semicontinuous. Thus the set
\[
C:=\{x\in \R: w(x)<e^x\lor K\}
\]
is open. From $(ii)$, we deduce that $C=D^c$ and therefore $D$ is a closed set.  The fact that  $w$ is non-decreasing and that $D$ is a closed set implies that there exists a $c^*\leq K$ such that $D=[c^*,\infty)$. In that case $\tau_D=\tau^+_{c^*}$.\\
$(iv)$ We first note that  from  the definition of $w$ as an infimum, we have
\[
\e\Big[U^{(y)}_{\tau^+_{c^*-x}}-U^{(y)}_{\tau^+_{c^*-y}}\Big]\ge 0.
\]
Now, for  $y\ge x$, it holds that $e^x\lor K-e^y\lor K\ge e^x-e^y$ and $\tau^+_{c^*-x}\ge \tau^+_{c^*-y}$. Therefore
\[
\begin{split}
w(x)-w(y)&=\e\Big[U^{(x)}_{\tau^+_{c^*-x}}-U^{(y)}_{\tau^+_{c^*-x}}+U^{(y)}_{\tau^+_{c^*-x}}-U^{(y)}_{\tau^+_{c^*-y}}\Big]\\
&\ge (e^x-e^y)\e\left[e^{-q\tau^+_{c^*-x}+X_{\tau^+_{c^*-x}}}+\beta\int_0^{\infty} e^{-qs+X_s}\ud s\right]+\e\Big[U^{(y)}_{\tau^+_{c^*-x}}-U^{(y)}_{\tau^+_{c^*-y}}\Big]\\
&\ge \mathfrak{K}(c^*,\beta)(e^x-e^y),
\end{split}
\]
for some constant $\mathfrak{K}(c^*,\beta)>0$ which depends on $c^*$ and $\beta$. Therefore, using part $(i)$, we deduce that $w$ is continuous and moreover that $w(c^*)=K$. \\
$(v)$ In what follows, for $q>0$, it is convenient to denote the function $w$ by $w(x,q)$ and $U_t=U_t(q)$ for all $t\ge 0$. Note that for any $t\ge 0$, $U_t(q)$ is non-increasing in $q$. Hence,
\[
w(x,q)\geq w(x,q_0)=\inf_{\sigma\in \mathcal{T}_{0,\infty}}\ex\Big[U_{\sigma\land \tau^+_{\log K}}(q_0)\Big]\qquad\textrm{for } q<q_0.
\]
On the other hand, recall from Theorem \ref{big-q} that when $q=q_0,$ a saddle point for the stochastic game (\ref{nash}) is given by $\tau^*=\sigma^*=\tau^+_{\log K}$, and in particular the value function satisfies
\[
V(x,q_0)=\ex\Big[U_{\tau^+_{\log K}}(q_0)\Big].
\]
Therefore, appealing to  the definition of $V$ as an infimum and using the lower bound on the solution to Theorem \ref{big-q}, we have
\[
w(x,q)\geq w(x,q_0)=M(\tau^+_{c^*}, \tau^+_{\log K})\ge V(x,q_0)\ge e^x.
\]
$(vi)$ and $(vii)$ These are trivial statements.\\
$(viii)$ and $(ix)$ These are standard results from the theory of optimal stopping. See for example Theorem 2.2 on p.29 or Theorem 2.4 p.37 of Peskir and Shiryaev \cite{PS}.
\hfill$\square$

\bigskip

According to the previous Lemma and the Verification Lemma, a stochastic saddle point of the K\"uhn--Gapeev game  exists and is given by $\tau^*=\tau^+_{\log K}$ and $\sigma^*=\tau^+_{c^*}$, for a given $c^*\le \log K$. (Note that the condition (\ref{condvl}) is dealt with in the same way as before). Therefore the associated value function  is given by
\[
V(x)=\ex\left[e^{-q\tau^+_{c^*}}\Big(e^{X_{ \tau^+_{c^*}}}\lor K\Big)+\int_0^{ \tau^+_{c^*}} e^{-qs}(\alpha+\beta e^{X_s})\ud s\right].
\]
The proof of Theorem \ref{medium-q} is thus complete as soon as we can characterize $c^*$ as given in the statement of the theorem.

 Suppose  that $b^2>0$. Our objective is to show that $\tau^*=\sigma^*=\tau^+_{\log K}$ is the stochastic saddle point provided $q$ is smaller than $q_0$ but not too small (to be made precise below). We again do this with the help of the Verification Lemma.

We show that $c^*=\log K$ if and only if $H_{\log K}'(\log K-)\geq 0.$
Note that from (\ref{prime}) we find that
\begin{equation}
H_{\log K}'(\log K - )=K+ \frac{2\alpha}{\Phi(q)b^2}\left(Ke^{-a^*(q)} -1\right),
\label{notsmooth}
\end{equation}
where we have used the fact that $W^{(q)}(0+)=0$ and $W^{(q)\prime}(0+) = 2/b^2$ when $b^2>0$ (cf. Chapter 8 of Kyprianou \cite{ky}).
Taking account of the monotonicity of $H_{\log K}(x, q)$ in $q$ this implies that those $q\in(0,q_0)$ for which
\[
K
+ \frac{2\alpha}{\Phi(q)b^2}\left(Ke^{-a^*(q)} -1  \right)\geq 0
\]
form an interval the left end point of which we shall denote by $q_1.$
First consider $q>q_1$. It then holds that $H_{\log K}'(\log K - )>0$ and hence $H_{\log K}(x,q)<H_{\log K}(\log K,q)=K$ for $x\in[\log K-\varepsilon,\log K)$ for some $\varepsilon>0$.
Now any choice of $c^*<\log K$ would imply $H_{\log K}(x,q)<w(x)$ for some $x<\log K,$ since $w(x)=K$ for all $x\in[c^*,\log K]$.
This leads to an immediate contradiction due to the fact that $\tau_{\log K}^+$ is a feasible strategy for the optimal stopping problem (\ref{opstop}). We conclude that for $q>q_1$ we have that $c^*=\log K$.

Next, we show that $c^*=\log K$ also in the case when $q=q_1$. For any $q>q_1$ it holds that $H_{\log K}(x,q)\leq K\vee e^x$ for all $x$ and thus we find that $H_{\log K}(x,q_1)\leq K\vee e^x$ due to continuity of $H_{\log K}(x,q)$ in $q$.

Furthermore, note that
\begin{equation}
\left(e^{-q(t\land \tau^+_{\log K})}w(X_{t\land \tau^+_{\log K}})+\int_0^{t\land \tau^+_{\log K}} e^{-qs}(\alpha+\beta e^{X_s})\ud s,\, t\ge 0\right)
\label{isamg}
\end{equation}
is a martingale for $q\in (q_1,q_0),$ as it is both a submartingale and a supermartingale due to items (viii) and (ix) of Lemma \ref{lemt3}. From monotone convergence it follows that (\ref{isamg}) is also a martingale when $q=q_1$.


Next, we show that $q_1>\alpha/K$. It seems unclear how to prove this inequality directly using the definition of $q_1$ and instead we argue by contradiction, hence suppose that $q_1\leq \alpha/K$. Due to monotonicity in $q$ in the definition of $w$ and Theorem \ref{qsmall} it would then follow that $K\vee e^x\geq w(x,q_1)\geq V(x,\alpha/K)=K\vee e^x$ for all $x$. Hence in this case
\[
\left(e^{-q_1(t\land \tau^+_{\log K})}(K\vee e^{X_{t\land \tau^+_{\log K}}})+\int_0^{t\land \tau^+_{\log K}} e^{-q_1s}(\alpha+\beta e^{X_s})\ud s,\, t\ge 0\right)
\]
is a martingale. Recall from the proof of Theorem \ref{qsmall} however  that when $x<\log K$, the process above is strictly increasing.
We get a contradiction with the martingale property and thus conclude that $q_1>\alpha/K$.

From (\ref{notsmooth}) it is clear that smooth pasting can only occur when $H_{\log K}'(\log K - ) = 0$ or $K$. This occurs precisely at the end points of the interval  $[q_1, q_0]$.

We conclude the proof by noting that when $b^2=0$, by considering (\ref{fit}) and (\ref{prime}) with $a=\log K$ and recalling that $W^{(q)}(0+)>0$ if $X$ has bounded variation and $W^{(q)\prime}(0+)=+\infty$ if $X$ has unbounded variation,  the  strategies $\tau^*=\sigma^* =\inf\{t\geq 0 : X_t >\log K\}$ do not constitute  a stochastic saddle point when $q<q_0$ as otherwise the necessary upper bound, $K\vee e^x$ on the value function $V$ will not be  respected.
\hfill$\square$

\section{Proof of Theorem \ref{small-q}}
The proof of Theorem \ref{small-q} again relies on the optimal stopping problem introduced in the previous Section. Assume that $\alpha/K<q<q_1$, which is possible thanks to Theorem \ref{medium-q}.

Let us first address the issue of continuous and smooth fit. We know from Lemma \ref{lemt3} that the value function $V$ is always continuous and hence in particular there is always continuous fit at the point $c^*$. Note that necessarily $c^*<\log K$ as otherwise $c^* = \log K$ and then from the previous theorem, $q=q_1$ which is a contradiction. As we shall see, this will be sufficient to uniquely characterize the value $c^*$ in the case that $X$ has paths of bounded variation. When $X$ has paths of unbounded variation, consistently with prior experience, continuous fit is not enough and the following lemma will be needed instead.

\begin{lemma}
When $X$ has paths of unbounded variation it holds that $V'(c^*-)= V'(c^*+)=0$.
\end{lemma}
\noindent {\it Proof:} Thanks to monotonicity of the value function we know that $V(x)\leq V(c^*)$ for all $x\leq c^*$ and hence
\[
\liminf_{x\uparrow c^*}\frac{V(c^*) - V(x)}{c^*-x}\geq 0.
\]
The proof is thus complete as soon as we show that
\begin{equation}
\limsup_{x\uparrow c^*}\frac{V(c^*) - V(x)}{c^*-x}\leq 0.
\label{limsup}
\end{equation}

To this end, let $\epsilon>0$ and introduce  $\tau^+_{c^*+\epsilon} = \inf\{t>0 : X_t > c^* + \epsilon \}$, $\tau^-_{c^*-\epsilon}= \inf\{t> 0 : X_t <c^*-\epsilon\}$ and $\tau = \tau^+_{c^*+\epsilon}\wedge \tau^-_{c^*-\epsilon}$. From parts $(iv)$ and $(viii)$ of  Lemma \ref{lemt3}, we have
\begin{eqnarray}
\lefteqn{\mathbb{E}_{c^*}\left[e^{-q\tau} V(X_\tau) + \int_0^{ \tau} e^{-qs}(\alpha+\beta e^{X_s})\ud s\right]}&&\notag\\
&&\geq V(c^*)\mathbb{E}_{c^*}\Big[e^{-q\tau^-_{c^*-\epsilon} } \mathbf{1}_{\{\tau^-_{c^*-\epsilon} <\tau^+_{c^*+\epsilon} \}}\Big] + K \Big(1- \mathbb{E}_{c^*}\Big[e^{-q\tau^-_{c^*-\epsilon}} \mathbf{1}_{\{\tau^-_{c^*-\epsilon}<\tau^+_{c^*+\epsilon}\}}\Big]\Big).\label{lower}
\end{eqnarray}
On the other hand, we have with the help of spectral positivity of $X$, $V(c^*)=K$ and the upper bound on $V$ that
\begin{eqnarray}
\lefteqn{
\mathbb{E}_{c^*}\left[e^{-q\tau} V(X_\tau) + \int_0^{ \tau} e^{-qs}(\alpha+\beta e^{X_s})\ud s\right]
}&&\notag\\
&\leq & V(c^*-\epsilon)\mathbb{E}_{c^*}\Big[e^{-q\tau^-_{c^*-\epsilon}} \mathbf{1}_{\{\tau^-_{c^*-\epsilon}<\tau^+_{c^*+\epsilon}\}}\Big]+ K\mathbb{E}_{c^*}\Big[e^{-q\tau^+_{c^*+\epsilon}} \mathbf{1}_{\{\tau^+_{\log K} \neq \tau^+_{c^*+\epsilon}<\tau^-_{c^*-\epsilon}\}}\Big]\notag\\
&&+ \mathbb{E}_{c^*}\Big[e^{-q\tau^+_{c^*+\epsilon}} e^{X_{\tau^+_{c^*+\epsilon}}}\mathbf{1}_{\{\tau^+_{\log K} = \tau^+_{c^*+\epsilon}<\tau^-_{c^*-\epsilon}\}}\Big]
+\mathbb{E}_{c^*}\left[\int_0^{\tau} e^{-qs}(\alpha+\beta e^{X_s})\ud s\right]\notag\\
&\leq& V(c^*-\epsilon)\mathbb{E}_{c^*}\Big[e^{-q\tau^-_{c^*-\epsilon}} \mathbf{1}_{\{\tau^-_{c^*-\epsilon}<\tau^+_{c^*+\epsilon}\}}\Big]+ K\mathbb{E}_{c^*}\Big[e^{-q\tau^+_{c^*+\epsilon}} \mathbf{1}_{\{\tau^+_{\log K} \neq \tau^+_{c^*+\epsilon}<\tau^-_{c^*-\epsilon}\}}\Big]\notag\\
&&+K\mathbb{P}_{c^*}^{(-1)}\Big(\tau^+_{\log K} =\tau^+_{c^*+\epsilon}<\tau^-_{c^*-\epsilon}\Big) + \frac{(\alpha + \beta e^{c^*+\epsilon})}{q}\Big(1 - \mathbb{E}_{c^*}\big[e^{-q\tau} \big]\Big)\label{upper}
\end{eqnarray}
Next, we claim that the last two terms on the right hand side above are $o(\epsilon)$. For the first of these two terms, the claim follows by Lemma 10 of Baurdoux and Kyprianou \cite{bky}. The second of these two terms is proportional to (cf. Chapter 8 of Kyprianou \cite{ky})
\[
1 - \mathbb{E}_{c^*}\big[e^{-q\tau} \big] = q\frac{W^{(q)}(\epsilon)}{W^{(q)}(2\epsilon)} \int_0^{2\epsilon} W^{(q)}(y){\rm d}y   - q  \int_0^{\epsilon} W^{(q)}(y){\rm d}y
\]
which is $o(\epsilon)$ on account of the fact that $W^{(q)}$ is monotone increasing with $W^{(q)}(0+)=0$ (the latter is due to the assumption that $X$ has paths of unbounded variation).

Taking this into account and combining the inequalities (\ref{lower}) and (\ref{upper}) we get
\begin{eqnarray*}
\frac{V(c^*)-V(c^*-\epsilon)}{\epsilon} &\leq& K \frac{ \mathbb{E}_{c^*}\big[e^{-q\tau} \big] -1}{\epsilon \mathbb{E}_{c^*}\Big[e^{-q\tau^-_{c^*-\epsilon}} \mathbf{1}_{\{\tau^-_{c^*-\epsilon}<\tau^+_{c^*+\epsilon}\}}\Big] } +\frac{o(\epsilon)}{\epsilon}\\
&=&\frac{qK}{\epsilon}  \left( \frac{W^{(q)}(2\epsilon)}{W^{(q)}(\epsilon)} \int_0^{\epsilon} W^{(q)}(y){\rm d}y
  -\int_0^{2\epsilon} W^{(q)}(y){\rm d}y\right) +\frac{o(\epsilon)}{\epsilon}.
      \end{eqnarray*}
Lemma 11 in Baurdoux and Kyprianou \cite{bky} states that $\limsup_{\epsilon\downarrow 0}W^{(q)}(2\epsilon)/W^{(q)}(\epsilon)\leq 2$ and hence the expression in the brackets on the right-hand side above is $o(\epsilon)$. This in turn implies (\ref{limsup}) and hence the proof is complete.\hfill$\square$

\bigskip

 Define, for each $c\leq \log K$, $G_c(x):=\mathbb{E}_x\big[U_{\tau^+_c}\big]$, that is to say
\begin{equation}\label{identu}
G_c(x) = \mathbb{E}_x\Big[e^{-q\tau^+_c}\big(e^{X_{\tau^+_c}}\lor K\big)\Big]+\int_0^{\infty}\mathbb{E}_x \Big[e^{-qs}(\alpha+\beta e^{X_s})\ind_{\{s\le \tau^+_c\}}\Big]\ud s,
\end{equation}
and note that for $x\ge c$, we have
\[
G_c(x)=e^x\lor K.
\]
We may now put the features of continuous and smooth fit to use and characterize the value of $c^*$.
Our immediate aim is  to give and explicit form of $G(x)$,  for $x<c$, in terms of  scale functions and the characteristics of $X$. We first note that the integral on the right-hand side of (\ref{identu}) has been computed before and  is equal to
\[
\int_{-\infty}^c \big(\alpha+\beta e^{y}\big)\Big(e^{-\Phi(q)(c-y)}W^{(q)}(c-x)-W^{(q)}(y-x)\Big)\ud y.
\]
The first term on the right-hand side of (\ref{identu}) satisfies
\begin{equation}\label{fuckeq}
\mathbb{E}_x\Big[e^{-q\tau^+_c}\big(e^{X_{\tau^+_c}}\lor K\big)\Big]=K\mathbb{E}_x\Big[e^{-q\tau^+_c }\Big]+\mathbb{E}_x\Big[e^{-q\tau^+_c}\big(e^{X_{\tau^+_c}}-K\big)\ind_{\{X_{\tau^+_{c}}>\log K\}}\Big].
\end{equation}
Recall that $\widehat{\p}$ denotes the law of $\widehat{X}=-X$. By Theorem 8.1 in \cite{ky}, we get that the first term on the right-hand side of (\ref{fuckeq}) satisfies
\[
\begin{split}
\mathbb{E}_x\Big[e^{-q\tau^+_c}\Big]&=\widehat{\mathbb{E}}_{c-x}\Big[e^{-q\tau^-_0}\Big]=Z^{(q)}(c-x)-\frac{q}{\Phi(q)}W^{(q)}(c-x).
\end{split}
\]
Now using the exponential change of measure (\ref{esscher}) with $\lambda =\Phi(q)$,  we write  the second term in the right-hand side of (\ref{fuckeq}) as follows
\[
\begin{split}
\mathbb{E}_x\Big[e^{-q\tau^+_c}\big(e^{X_{\tau^+_c}}-K\big)&\ind_{\{X_{\tau^+_{c}}>\log K\}}\Big]=\mathbb{E}_x^{\Phi(q)}\Big[e^{\Phi(q)(X_{\tau^+_c}-x)}\big(e^{X_{\tau^+_c}}-K\big)\ind_{\{X_{\tau^+_{c}}>\log K\}}\Big]\\
&=\mathbb{E}^{\Phi(q)}\Big[e^{\Phi(q)X_{\tau^+_{c-x}}}\big(e^{x+X_{\tau^+_{c-x}}}-K\big)\ind_{\{X_{\tau^+_{c-x}}+x>\log K\}}\Big].
\end{split}
\]
Let $f(y)=e^{\Phi(q)y}\big(e^{x+y}-K\big)\ind_{\{y+x>\log K\}}.$
From Theorem 4.4 in \cite{ky} and since $x<c \le \log K$, we deduce
\[
\begin{split}
\mathbb{E}^{\Phi(q)}&\Big[f\big(X_{\tau^+_{c-x}}\big)\Big]=\mathbb{E}^{\Phi(q)}\Big[f\big(X_{\tau^+_{c-x}}\big)\ind_{\{X_{\tau^+_{c-x}}>c-x\}}\Big]\\
&=\mathbb{E}^{\Phi(q)}\left[\int_0^\infty \ud t f\big(X_{t}\big)\ind_{\{\overline{X}_{t-}<c-x\}}\ind_{\{X_{t}>c-x\}}\right]\\
&=\mathbb{E}^{\Phi(q)}_{x-c}\left[\int_0^\infty\ud t\int_0^\infty\Pi^{\Phi(q)}(\ud u)  f\big(u+X_{t-}+c-x\big)\ind_{\{\tau^+_0>t\}}\ind_{\{u+X_{t-}>0\}}\right]\\
&=\int_0^\infty\Pi^{\Phi(q)}(\ud u) \int_0^\infty\ud t\,\mathbb{E}^{\Phi(q)}_{x-c}\left[ f\big(u+X_{t-}+c-x\big)\ind_{\{\tau^+_0>t\}}\ind_{\{u+X_{t-}>0\}}\right]\\
&=\int_0^\infty\Pi^{\Phi(q)}(\ud u) \int_{-\infty}^0 R_{\Phi(q)}(x-c,\ud y; 0)\ f\big(u+y+c-x\big)\ind_{\{u+y>0\}}
\end{split}
\]
where
$R_{\Phi(q)}(z,\ud y; 0) $ plays the role of $R(z,\ud y;0) $ but under the measure $\mathbb{P}^{\Phi(q)}_z$. Therefore, by Corollary 8.8 in Kyprianou \cite{ky} we get
\[
\begin{split}
&\mathbb{E}^{\Phi(q)}\Big[f\big(X_{\tau^+_{c-x}}\big)\Big]\\
&=\int_0^\infty\Pi^{\Phi(q)}(\ud u) \int_{-\infty}^0 \ud y \Big(W_{\Phi(q)}(c-x)-W_{\Phi(q)}(c-x+y)\Big) f\big(u+y+c-x\big)\ind_{\{u+y>0\}}.
\end{split}
\]
Finally, putting the pieces together, using in particular that  $\Pi^{\Phi(q)}({\rm d} x) = e^{-\Phi(q)x}\Pi({\rm d}x)$ and $W_{\Phi(q)}(x) = e^{-\Phi(q)x}W^{(q)}(x)$, we obtain the following formula for $G_c(x)$, when $x<c$,
\[
\begin{split}
G_c(x)&=K\Big(Z^{(q)}(c-x)-\frac{q}{\Phi(q)}W^{(q)}(c-x)\Big)\\
&+\int_{-\infty}^c \big(\alpha+\beta e^{y}\big)\Big(e^{-\Phi(q)(c-y)}W^{(q)}(c-x)-W^{(q)}(y-x)\Big)\ud y\\
&+e^{\Phi(q)(x-c)}\int_0^\infty e^{-\Phi(q)u}\Pi(\ud u) \int_{-\infty}^0 \ud y \Big(W^{(q)}(c-x)-e^{-\Phi(q)y}W^{(q)}(c-x+y)\Big)\\
&\hspace{5.3cm}\times e^{\Phi(q)(u+y)}\big(e^{c+u+y}-K\big)\ind_{\{u+y+c>\log K\}}.
\end{split}
\]

Now that we have an expression for $G_c$ we may find the one which corresponds to the optimal solution by choosing $c= c^*$ so that there is smooth or continuous fit accordingly with the path variation of $X$.

\bigskip

{\bf Bounded variation case:} In this case we know that $W^{(q)}(0+)=1/d>0$. Hence, checking for a discontinuity at $c$ we find that
\begin{equation}\label{Gc-}
\begin{split}
G_c(c-)&=K\left(1-\frac{q}{\Phi(q)}\frac{1}{d}\right)+\frac{e^{-\Phi(q)c}}{d}\left(\frac{\alpha}{\Phi(q)}e^{\Phi(q)c}+\frac{\beta}{\Phi(q)+1}e^{(\Phi(q)+1)c}\right)\\
&\qquad +\frac{K}{d}\int_{0}^\infty \Pi(\ud z+\log K-c)f(z),
\end{split}
\end{equation}
where
\begin{equation}\label{intf}
f(z)=\left \{
 \begin{array}{ll}
\frac{1}{\Phi(q)+1}e^z(1-e^{-(\Phi(q)+1)z})-\frac{1}{\Phi(q)}(1-e^{-\Phi(q)z})&\textrm{if $z\ge 0$}\\
0 & \textrm{if $z<0$.}
 \end{array} \right.
\end{equation}
It is important to note that
\[
f(z)\sim z^2\qquad\textrm{as $z\to 0$, }\quad\textrm{ and } \quad f(z)\sim\frac{1}{\Phi(q)+1}e^z\qquad\textrm{as $z\to \infty$},
\]
thus from the hypothesis (A) and the fact that  $\Pi$ is a L\'evy measure, we have
\[
\int_{0}^\infty \Pi(\ud z)f(z)<\infty.
\]
So, we take
\begin{equation}\label{inteq1}
\begin{split}
\frac{q}{\Phi(q)}&=\frac{1}{K}\left(\frac{\alpha}{\Phi(q)}+\frac{\beta}{\Phi(q)+1}e^{c}\right)\\
&\qquad-\frac{1}{\Phi(q)}\int_{0}^\infty \Pi(\ud z+\log K-c)(1-e^{-\Phi(q)z})\\
&\qquad+\frac{1}{\Phi(q)+1}\int_{0}^\infty \Pi(\ud z+\log K-c)e^z(1-e^{-(\Phi(q)+1)z}).
\end{split}
\end{equation}
In order to show that this expression has a unique solution, it is more convenient to note from (\ref{Gc-}) that
\[
\lim_{c\downarrow - \infty}G_c(c-) = K-\frac{Kq}{\Phi(q)d}\left(1-\frac{\alpha}{Kq}\right)<K
\]
on account of the assumption that $q>\alpha/K$. Moreover, as in the case of bounded variation paths,
\[
\psi(\theta) = d\theta - \int_0^\infty (1-e^{-\theta x})\Pi({\rm d}x),
\]
and $\psi(\Phi(q)) = q$, we may compute from (\ref{Gc-})
\begin{eqnarray*}
\lim_{c\uparrow\log K}G_c(c-)&=&K+\frac{1}{d }\left(\frac{\alpha}{\Phi(q)} -K\frac{(q-\psi(-1) -\beta )}{(\Phi(q)+1)}  \right)>K,
\end{eqnarray*}
where the strict inequality follows from the fact  $q<q_1=q_0$. Thus, we get the existence of the unique solution if we prove  that $G_\cdot(\cdot -)$ is continuous and increasing in $(-\infty,\log K]$.  The continuity of $G_.(.-)$ follows from (\ref{Gc-}) and the fact that when the measure $\Pi$ has an atom at $\log K-c$, the integrand on the right-hand side of (\ref{Gc-}) is equal to $0$ at $z=0$.

Now, note that
\[
f'(z)=\frac{1}{\Phi(q)+1}(e^z-e^{-\Phi(q)z})> 0\qquad \textrm{for all } z> 0,
\]
which implies that $f$ is positive and increasing. Then from (\ref{Gc-}) it is clear that  $G_.(.-)$ is increasing in $(-\infty,\log K]$.

\bigskip
{\bf Unbounded variation case:} In this case $W^{(q)}(0+)=0$ and hence in the above analysis one sees that $G_c(c-)=K = G_c(c+)$.  In that case, the principle of smooth fit can be implemented and we insist on choosing $c$ such that there is no discontinuity in $G_c'(c-)$.

Recall that $W^{(p)}\in C^1(0,\infty)$ and let $x<c$. Therefore, using a standard argument involving dominated convergence to  differentiate through the integral in the last term of $G_c(x)$, we have
\[
\begin{split}
G_c'(x)&= K\left(\frac{q}{\Phi(q)}W^{(q)\prime}(c-x)-qW^{(q)}(c-x)\right)\\
&\quad+\int_{-\infty}^c \big(\alpha+\beta e^{y}\big)\big(W^{(q)\prime}(y-x)-W^{(q)\prime}(c-x) e^{-\Phi(q)(c-y)}\big)\ud y \\
&\qquad+\Phi(q)e^{\Phi(q)(x-c)}\int_0^\infty e^{-\Phi(q)u}\Pi(\ud u) \int_{-\infty}^0 \ud y \Big(W^{(q)}(c-x)-e^{-\Phi(q)y}W^{(q)}(c-x+y)\Big)\\
&\qquad\quad\times e^{\Phi(q)(u+y)}\big(e^{c+u+y}-K\big)\ind_{\{u+y+c>\log K\}}\\
&\qquad\qquad+e^{\Phi(q)(x-c)}\int_0^\infty e^{-\Phi(q)u}\Pi(\ud u) \int_{-\infty}^0 \ud y \Big(e^{-\Phi(q)y}W^{(q)\prime}(c-x+y)-W^{(q)\prime}(c-x)\Big)\\
&\qquad\qquad\quad\times e^{\Phi(q)(u+y)}\big(e^{c+u+y}-K\big)\ind_{\{u+y+c>\log K\}}.
\end{split}
\]
Also, recall that $W^{(q)}(0+)=0$ and that $W^{(q)\prime}(0+)=2/b^2$ which should be interpreted as $+\infty$ in the case that the Gaussian coefficient $b^2 =0$,
\[
\begin{split}
G'_c(c-)&=K\frac{q}{\Phi(q)}W^{(q)\prime}(0+)-W^{(q)\prime}(0+) e^{-\Phi(q)c}\int_{-\infty}^c \big(\alpha+\beta e^{y}\big)e^{\Phi(q)y}\ud y\\
&\quad-W^{(q)\prime}(0+)\int_0^\infty e^{-\Phi(q)u}\Pi(\ud u) \int_{-\infty}^0 \ud y e^{\Phi(q)(u+y)}\\
&\qquad\times \big(e^{c+u+y}-K\big)\ind_{\{u+y+c>\log K\}}.
\end{split}
\]
In order to obtain the smooth fit $G'_c(c+)=0$ we necessarily must have:
\[
\begin{split}
\frac{q}{\Phi(q)}&= \frac{e^{-\Phi(q)c}}{K}\int_{-\infty}^c \big(\alpha+\beta e^{y}\big)e^{\Phi(q)y}\ud y\\
&\quad+\frac{1}K\int_0^\infty e^{-\Phi(q)u}\Pi(\ud u) \int_{-\infty}^0 \ud y e^{\Phi(q)(u+y)}\\
&\qquad\times \big(e^{c+u+y}-K\big)\ind_{\{u+y+c>\log K\}}.
\end{split}
\]
After some algebra, we get
\begin{eqnarray}
\frac{q}{\Phi(q)}&=&\frac{1}{K}\left(\frac{\alpha}{\Phi(q)}+\frac{\beta}{\Phi(q)+1}e^{c}\right)\nonumber\\
&&-\frac{1}{\Phi(q)}\int_{0}^\infty \Pi(\ud z+\log K-c)(1-e^{-\Phi(q)z})\nonumber\\
&&\quad+\frac{1}{\Phi(q)+1}\int_{0}^\infty \Pi(\ud z+\log K-c)e^z(1-e^{-(\Phi(q)+1)z})\label{equationuv},
\end{eqnarray}
which is the same identity as in (\ref{inteq1}). In order to prove the existence of a unique solution of the above identity we will follow similar arguments as those used in the bounded variation case. Let us define
\begin{equation}\label{Fc}
\begin{split}
F(c)&=K\left(1-\frac{q}{\Phi(q)}\right)+e^{-\Phi(q)c}\left(\frac{\alpha}{\Phi(q)}e^{\Phi(q)c}+\frac{\beta}{\Phi(q)+1}e^{(\Phi(q)+1)c}\right)\\
&\qquad +K\int_{0}^\infty \Pi(\ud z+\log K-c)f(z),
\end{split}
\end{equation}
where $f$ is defined as in  (\ref{intf}). Note that $c$ is a solution to (\ref{equationuv}) if and only if $c$ solves $F(c)=K$. Similarly to the bounded variation case, we have that
\[
\int_{0}^\infty \Pi(\ud z)f(z)<\infty.
\]
Now, we note from (\ref{Fc}) that
\[
\lim_{c\downarrow - \infty}F(c) = K-\frac{Kq}{\Phi(q)}\left(1-\frac{\alpha}{Kq}\right)<K
\]
on account of the assumption that $q>\alpha/K$. Moreover, recall that
\[
\psi(\theta) = a\theta +\frac{b^2}{2}\theta^2+\int_0^\infty \big(e^{-\theta x}-1+\theta x\ind_{\{x<1\}}\big)\Pi({\rm d}x),
\]
and $\psi(\Phi(q)) = q$, then after some straightforward computations we get
\[
\int_{0}^\infty \Pi(\ud z)f(z)=\frac{\psi(-1)+a-b^2/2}{\Phi(q)+1}+\frac{q-a\Phi(q)-b^2/2\Phi^2(q)}{\Phi(q)(\Phi(q)+1)}.
\]
Hence from (\ref{Fc})
\begin{eqnarray*}
\lim_{c\uparrow\log K}F(c)&=&K+\left(\frac{\alpha}{\Phi(q)} -K\frac{(q-\psi(-1) -\beta )}{(\Phi(q)+1)}-K\frac{ b^2}{2}  \right)>K,
\end{eqnarray*}
where the strict inequality follows from the fact  $q<q_1$ (recall that $q_1=q_0$ when $b^2=0$). The existence of the unique solution now follows from the continuity and the monotonicity of $F$ which can be proved as in the bounded variation case.
\section*{Acknowledgments}
This research was funded by EPSRC grant EP/D045460/1. We are grateful to Kees van Schaik for useful discussions.

\end{document}